\documentclass[11pt]{amsart}

\theoremstyle{definition}

\theoremstyle{remark}

\numberwithin{equation}{section}

\begin{document}

\title[On a Question of Bouras ]
 {On a Question of Bouras concerning weak compactness of almost Dunford-Pettis sets}
\author[J. X. Chen]
{Jin Xi Chen }

\address{Department of Mathematics, Southwest Jiaotong
University, Chengdu 610031, PR China}
 \email{jinxichen@home.swjtu.edu.cn}

\author[L. Li]
{Lei Li}
\address{Department of Mathematics, Nankai University, Tianjin 300071, PR China}
\email{leilee@nankai.edu.cn}

\thanks{ This work was supported  by NSFC (No.11301285). The first author  was also supported by China Scholarship Council (CSC) and was visiting the University of South Carolina when this work was completed.}

%    General info
%\dedicatory{}


\begin{abstract}
We give a positive answer to the question of K. Bouras [`Almost Dunford-Pettis sets in Banach lattices', \textit{Rend. Circ. Mat. Palermo (2)} \textbf{ 62} (2013), 227--236] concerning weak compactness of almost Dunford-Pettis sets in Banach lattices. That is, every almost Dunford-Pettis set in a Banach lattice $E$ is relatively weakly compact  if, and only if, $E$ is a $KB$-space.
\end{abstract}

\maketitle \baselineskip 5.3mm
\par  Let $E$ be a Banach lattice. Recall that a norm bounded subset $A$ of $E$ is called to be \textit{$L\,$-weakly compact} if $\|x_n\|\rightarrow0$ for every disjoint sequence $(x_{n})$ contained in the solid hull of $A$. Every $L\,$-weakly compact set is relatively weakly compact set, but the converse does not hold in general.  In an $L$-space, $L\,$-weakly compact sets and relatively weakly compact sets coincide. More generally, every relatively weakly compact subset of $E$ is $L$-weakly compact if, and only if, $E$ has the positive Schur property. Hereby, we say a Banach lattice $E$ has the \textit{positive Schur property} if every weakly null sequence with positive terms in $E$ is norm null.
 \par Following Bouras \cite{Kh}, a bounded subset $A$ of the Banach lattice $E$ is said to be an \textit{almost Dunford-Pettis set} if every disjoint weakly null sequence $(f_{n})$ of $ E^{\,\prime}$ converges uniformly to zero on $A$,  that is, $\sup_{x\in A}|f_{n}(x)|\rightarrow0$. In \cite{Kh}, Bouras  showed that every $L\,$-weakly compact set in a Banach lattice $E$ is necessarily almost Dunford-Pettis. Also,  every relatively weakly compact set in $E$ is almost Dunford-Pettis if and only if  $E$ has the weak Dunford-Pettis property. If every almost Dunford-Pettis set in $E$ is relatively weakly compact (in particular, $L$-weakly compact), then $E$ is a $KB$-space \cite[Theorem 2.10]{Kh}. Conversely, if  $E$ is an $L$-space, or a dual $KB$-space, or the norm of $E^{\prime\prime}$ is order continuous, then  every almost Dunford-Pettis set in $E$ is $L$-weakly compact (and hence relatively weakly compact). In those three cases, $E$ is a $KB$-space. Therefore, Bouras posed the following open question.

\noindent \textbf{Question} \cite{Kh}. Does the assumption ``$E$ is a $KB$-space" imply that each almost Dunford-Pettis set in $E$ is relatively weakly compact (resp. $L$-weakly compact)?

In this short note, we give a positive answer to the question with respect to weak compactness of almost Dunford-Pettis sets. That is, every almost Dunford-Pettis set in a Banach lattice $E$ is relatively weakly compact  if, and only if, $E$ is a $KB$-space.
\par  For Banach lattice theory, we refer  the reader to \cite{AB, M} and also the original paper of Bouras \cite{Kh}.

\noindent \textbf{Theorem.} \textit{Let $E$ be a Banach lattice. Every almost Dunford-Pettis set in $E$ is relatively weakly compact  if, and only if, $E$ is a $KB$-space.}

\begin{proof} For the proof of the ``only if\," part see Theorem 2.10 of \cite{Kh}. We  need only to prove the ``if\," part. To this end, let $E$ be a $KB$-space and $A$ an almost Dunford-Pettis set in $E$.
\par First, we claim that the solid hull $sol(A)$ of $A$  is likewise almost Dunford-Pettis. Otherwise, there would exist a disjoint weakly null sequence $(f_{n})\subseteq E^{\,\prime}$ such that $\sup_{x\in sol(A)}|f_{n}(x)|>\varepsilon_0$ for some $\varepsilon_0>0$ and all $n\in \mathbb{N}$. So, we can find two sequences $(x_n)\subseteq A$ and $(y_n)\subseteq sol(A)$ satisfying$$|y_n|\,\leq\mid x_n|,\qquad 0<\varepsilon_0<| f_{n}(y_{n})|\leq |f_{n}|(|x_{n}|)$$for each $n$. Again, by \cite[Theorem 1.23]{AB} there exists a sequence $(g_n)\subseteq E^{\,\prime}$ such that $$|g_{n}|\leq|f_{n}|,\qquad |g_{n}(x_{n})|>\varepsilon_{0}.$$Clearly, $(g_{n})$ is a disjoint sequence. From the weak convergence of $(f_{n})$ it follows that $g_{n}\xrightarrow {w}0$ (cf. \cite[Thoerem 4.34]{AB}). Since $A$ is almost Dunford-Pettis, we have $$\varepsilon_{0}<|g_{n}(x_{n})|\leq\sup_{x\in A}|g_{n}(x)|\rightarrow0$$which is impossible. This proves that $sol(A)$ is an almost Dunford-Pettis set.
\par Now we can assume without loss of generality that $A$ is a solid almost Dunford-Pettis set in $E$. Let $\rho_{A}(f):=\sup_{x\in A}\langle|f|,\,|\,x|\rangle$ for each $f\in E^{\,\prime}$. Clearly, $\rho_{A}$ is a lattice seminorm on $E^{\,\prime}$, and $\rho_{A}(f)=\sup_{x\in A}|f(x)|$. Let $(f_n)$ be an arbitrary order bounded disjoint sequence in $E^{\,\prime}$. Since $f_n\xrightarrow {w}0$ and $A$ is almost Dunford-Pettis, we have $\rho_{A}(f_n) \rightarrow0$. Therefore, for every $0\leq f\in E^{\,\prime}$, by \cite[Theorem 2.3.3]{M} $A$ is approximated order bounded with respect to $\rho_{f}$, that is to say, for every $\varepsilon>0$ there exists $0\leq x\in E$ satisfying $$A\subseteq[-x,\,x]+\varepsilon B_{\rho_{f}}\eqno{(\ast)}$$where $B_{\rho_{f}}=\{x\in E:\rho_{f}(x)=f(|\,x|)\leq1\}$.
\par Note that the $w^{\ast}$-closure  $\overline{A}^{\,w^*}$ of $A$ in $E^{\prime\prime}$ is $w^{*}$-compact. To prove that $A$ is relatively weakly compact in $E$ it suffices to show that $\overline{A}^{\,w^*}\subseteq E$. Since $E$ is a $KB$-space, we have $E=(E^{\,\prime})^{\,\prime}_{n}$ (cf. \cite[Theorem 4.60]{AB}). Thus, we have to show that every element $x^{\,\prime\prime}$ of $\overline{A}^{\,w^*} $ is an order continuous linear functional on $E^{\,\prime}$, equivalently, $|\,x^{\,\prime\prime}|\in(E^{\,\prime})^{\,\prime}_{n} $. To this end, let $(f_{\alpha})$ be an arbitrary decreasing net in $(E^{\,\prime})^+$ with $f_{\alpha}\downarrow0$. It is enough to show that $\langle|\,x^{\,\prime\prime}|,\,f_{\alpha}\rangle\downarrow0$. For each $\alpha$ we have
\begin{eqnarray*}\langle|\,x^{\,\prime\prime}|,\,f_{\alpha}\rangle=\sup\{|\langle x^{\,\prime\prime},\,g\rangle|:|\,g|\leq f_{\alpha}\}&\leq& \sup\{|\langle g,\,x\rangle|:|\,g|\leq f_{\alpha},\,x\in A\}\\&\leq&\rho_{A}(f_{\alpha})
\end{eqnarray*}
 To finish the proof, we need only to prove that $\rho_{A}(f_{\alpha})\downarrow0$. Indeed, it may be assumed that there exists an element $f\in E^{\,\prime}$ such that $0\leq f_{\alpha}\leq f $ for all $\alpha$. Let $\varepsilon>0$ be fixed. By $(\ast)$ we have $\rho_{A}(f_{\alpha})\leq f_{\alpha}(x)+\varepsilon$ for all $\alpha$. Note that $f_{\alpha}(x)\downarrow0$. It follows that $\inf_{\alpha}\rho_{A}(f_{\alpha})\leq\varepsilon$. Since $\varepsilon>0$ is abitrary, we have $\inf_{\alpha}\rho_{A}(f_{\alpha})=0$, as desired.
\end{proof}
\par In a $KB$-space without the weak Dunford-Pettis property (e.g., $\ell_{p}$, $1<p<\infty$) almost Dunford-Pettis sets and relatively weakly compact sets can not coincide.
\par There is another topic included in Bouras' question:
\par (Q1)\quad\textit{Is every almost Dunford-Pettis set in a KB-space L-weakly compact?}

\noindent We are not able to answer it and have to leave it still open. However, we can make some comments on (Q1). In \cite{CCJ} Chen \textit{et al.} introduced the class of almost limited sets in Banach lattices. A norm bounded subset $A$ of a Banach lattice $E$ is said to be an \textit{almost limited set} if every disjoint, weak$^{*}$ null sequence $(f_{n})$ of $ E^{\,\prime}$ converges uniformly to zero on $A$.  Clearly, every almost limited set in $E$ is an almost Dunford-Pettis set. Almost limited sets  and  $L$-weakly compact sets coincide in $E$ if and only if the norm of $E$ is order continuous (\cite[Theorem 2.6(2)]{CCJ}). Therefore, (Q1) is the same as the following question:
\par (Q1$^{\,\prime}$) \textit{Do almost Dunford-Pettis sets and almost limited sets   coincide in a KB-space?}
\par Let us recall that a Banach lattice $E$ has the positive Schur property   if, and only if, every relatively weakly compact subset of $E$ is $L$-weakly compact. We know that every Banach lattice with the positive Schur property is a $KB$-space with the weak Dunford-Pettis property. In 1994 Wnuk \cite {W4} left a  question which is still open until now:
\par (Q2)\quad \textit{Does every $KB$-space with the weak Dunford-Pettis property have the positive Schur property?}

\noindent By the Theorem of the present paper,  in a $KB$-space with the weak Dunford-Pettis property  almost Dunford-Pettis sets and relatively weakly compact sets  coincide. Hence, if the answer to (Q1) is positive, then  the answer to (Q2) is likewise positive. In other words, if the answer to (Q2) were  negative, then there would exist a $KB$-space where we can find an almost Dunford-Pettis set which is not $L$-weakly compact.

\end{document}